\documentclass[a4paper,twoside,10pt]{amsart}
\usepackage{amsmath,amsfonts,amssymb,amsthm,mathtools}
\allowdisplaybreaks
\usepackage{todonotes}
\newtheorem{theorem}{Theorem}[section]
\newtheorem{lemma}{Lemma}[section]

\usepackage{color,graphicx, float}
\usepackage[colorlinks=true,pdfstartview=FitH,pdfview=FitH]{hyperref}
\numberwithin{equation}{section}

\newcommand\mytop[2]{\genfrac{}{}{0pt}{}{#1}{#2}}
\newcommand{\secondsol}{\K{\expe^{\pi\iunit}q^2z}}

\def\expe{{\rm e}}
\def\iunit{{\rm i}}

\def\qpr#1#2#3{\left(#1;#2\right)_{#3}}
\def\sfrac#1#2{#1/#2}
\def\Eq#1{{\rm E}_q\left(#1\right)}
\def\bigO{{\mathcal O}}

\def\qphi#1#2#3#4#5#6{{{}_{#1}\phi_{#2}}\!\left(\genfrac{}{}{0pt}{}{#3}{#4};#5,#6\right)}

\def\qpsi#1#2#3#4#5#6{{{}_{#1}\psi_{#2}}\!\left(\genfrac{}{}{0pt}{}{#3}{#4};#5,#6\right)}
\def\qhat{{\hat{q}}}
\def\Thetaq#1{{\rm \theta}_q\left(#1\right)}
\def\dThetaq#1{{\rm \theta}'_q\left(#1\right)}
\def\Pq#1{\mathcal{P}_q\left(#1\right)}
\def\dPq#1{\mathcal{P}'_q\left(#1\right)}
\def\Cq{{C_q}}

\def\UqK#1#2#3#4{{\rm U}_{q,#1}\left(#2,#3;#4\right)}
\def\intd{{\rm d}}
\def\LqK#1{\mathcal{L}_{q,#1}}

\def\K#1{\kappa\left(#1\right)}

\def\ee{e}

\title{Second-Order Linear $q$-Difference Equations}
\author{Daniel Meikle}
\author{Adri Olde Daalhuis}
\address{School of Mathematics and Maxwell Institute for Mathematical Sciences, University of Edinburgh, Peter Guthrie Tait Road, Edinburgh, EH9 3FD, UK}
\email{D.Meikle@ed.ac.uk}
\email{A.OldeDaalhuis@ed.ac.uk}
\date{\today}

\begin{document}

\maketitle

\begin{abstract}
    Second–order $q$-difference equations arise in various areas such as basic hypergeometric series, orthogonal polynomials and the analytic theory of special functions. When the independent variable is scaled by powers of $q$ ($0<q<1$), such equations display asymptotic behaviour that differs markedly from their differential or classical difference counterparts. Our purpose is to give a complete, self-contained treatment of linear, second-order $q$-difference equations with divergent formal solutions, including explicit error bounds, exact resummations and Stokes phenomena.
\end{abstract}

\section{Introduction}\label{S:intro}

Second-order $q$-difference equations with an irregular singularity arise throughout
the theory of basic hypergeometric series, orthogonal polynomials and special
functions. For existence results for large classes of $q$-difference equations, we refer
to \cite{Ramis2013}. Here, we concentrate on three-term, linear $q$-difference
equations possessing \emph{divergent} formal solutions: the $q$-analogues of
second-order linear differential equations with an irregular singularity of rank one.
Our aim is a complete, self-contained treatment with explicit error bounds, exact
resummations, and a detailed description of the associated Stokes phenomena.

\subsection{Notation} We take $0<q<1$, $\qhat=\expe^{\frac{2\pi^2}{\ln q}}$ and 
define the $q$-Pochhammer symbol as $\qpr{a}{q}{\infty}=\prod\limits_{n=0}^\infty (1-aq^n)$ 
with $\qpr{a}{q}{\nu}=\qpr{a}{q}{\infty}/\qpr{aq^\nu}{q}{\infty}$, $a,\nu\in\mathbb{C}$. 
For a positive integer $n\in\mathbb{N}$, we therefore have 
$\qpr{a}{q}{n}=(1-a)(1-aq)\ldots(1-aq^{n-1})$ and 
$\qpr{a}{q}{-n}=\left(-q/a\right)^nq^{\binom{n}{2}}/\qpr{q/a}{q}{n}$. 
In our $q$-Laplace transforms, the alternatives for the exponential function are
\begin{equation}\label{Eq}
    \Eq{t}=\exp\left(\frac{\ln^2(t/\sqrt{q})}{2\ln q}\right)=q^{\frac{1}{8}}t^{-\frac{1}{2}}\exp\left(\frac{\ln^2t}{2\ln q}\right),
\end{equation}
and the reciprocal of the theta function 
\begin{equation}\label{Thetaq}
    \Thetaq{t}=\qpr{q,-t,-q/t}{q}{\infty}=\sum_{n=-\infty}^{\infty}q^{\binom{n}{2}}t^n.
\end{equation}
Since $\Eq{q^x}=q^{\frac12\left(x-\frac12\right)^2}$, terms of size $\Eq{t}$ are
called \emph{$q$-exponentially small}. 
Note that $\Eq{t}$ is $q$-exponentially small on all Riemann sheets!
We denote
\begin{equation}\label{X}
\begin{split}
    &\K{t}=\Cq \Eq{t}\qquad{\rm if}\quad \ee=E,\\
    &\K{t}=\frac{-1}{\ln(q)\Thetaq{t}}\,\quad{\rm if}\quad \ee=\theta,
\end{split}
\end{equation}
with $\Cq=1/\sqrt{-2\pi\ln q}$. For both choices, we have
\begin{equation} \label{Fprop}
    \K{t^{-1}}=t \K{t}, \qquad\K{q^nt}=t^n q^{\binom{n}{2}}\K{t},\qquad
    \int_0^\infty t^{n-1}\K{t}\intd t=q^{-\binom{n}{2}},
\end{equation}
with constraints, $n\in\mathbb{Z}$ if $\ee=\theta$, and $n\in\mathbb{C}$ if $\ee=E$. 

A formal series $\hat{w}(z)=\sum\limits_{n=0}^\infty w_nz^n$ is of \emph{$q$-Gevrey order}
$k\in{\mathbb R}^*$ if $|w_n|<q^{-k\binom{n}{2}}/F^n$ for some $F>0$; see \cite{Ramis:1992:GEF} for the growth of solutions of $q$-difference equations.

\subsection{The equation and its formal solution}
Our main $q$-difference equation is
\begin{equation} \label{newexpandedequation}
zy(zq^{-2})+a(z) y(zq^{-1})-b(z) y(z)=0,
\end{equation}
where both
\begin{equation}\label{abz}
a(z)=\sum_{m=0}^\infty a_mz^m,\qquad
b(z)=\sum_{m=0}^\infty b_mz^m,
\end{equation}
have a non-zero radius of convergence and $a_0b_0\not=0$. In \S\ref{S:general}, we
show that one may assume $a_0=b_0=1$ without loss of generality. The formal series
solution is
\begin{equation} \label{originaly}
    \hat{y}(z)=\sum_{n=0}^\infty c_n  q^{-\binom{n}{2}} \left(-qz\right)^n,
\end{equation}
with $c_0=1$ and, for $n\ge1$, the recurrence
\begin{equation} \label{cnrecurrence}
    c_n=\frac{1}{1-q^n}\left(c_{n-1}-\sum_{m=0}^{n-1}c_m\left(-1\right)^{n-m}q^{\binom{n}{2}-\binom{m}{2}}
    (a_{n-m}-b_{n-m}q^m)\right).
\end{equation}
If $C>0$ satisfies $|a_n|+|b_n|\leq C^n$ for $n\geq1$, then (Appendix~\ref{AppLemma})
$|c_n|\leq\qpr{-2C}{q}{n}/\qpr{q}{q}{n}$. In \S\ref{S:CH}, we will show that $\lim_{n\to\infty}c_n=K(q)$, the Stokes multiplier. The factor $q^{-\binom{n}{2}}$ in
\eqref{originaly} dominates, so the series diverges and is of $q$-Gevrey order $1$. 
Compared with a classical Gevrey-$1$
series, the divergence is much stronger, necessitating the use of different
Borel and Laplace transforms. 
The $q$-Borel transform of the right-hand side of \eqref{originaly} is
\begin{equation} \label{qBorel}
    B_q(t)=\sum_{n=0}^\infty c_n \left(-qt\right)^n.
\end{equation}
In \S\ref{S:vn}, we will show that this series converges for $|tq|<1$ and that $B_q(t)$ is a meromorphic
function. Its poles are simple poles at $t=-q^{-1},-q^{-2},-q^{-3},\ldots$,
and when we bound $t$ away from these poles,
\begin{equation}\label{Bqest}
    B_q(t)=\bigO\left(\left|t\right|^{3/2}\exp\left(\frac{\ln^2(D|t|)}{-2\ln q}  \right)  \right),
    \qquad \text{as}\quad |t|\to\infty,
\end{equation}
with $D$ a positive constant. This growth condition is needed for the $q$-Laplace transforms
in \S\ref{S:prelim} to converge.

\subsection{Main results}
Our principal contributions are the following:

\emph{(i) A convergent second solution.} In \S\ref{S:secondsolution}, we show that a
second solution of \eqref{newexpandedequation} has the form
$y(z)=\secondsol\,u(z)$, where $u(z)$ has a \emph{convergent} Taylor expansion about
the origin. This is a remarkable difference from the classical theory, where for a
rank-one irregular singularity, all formal solutions typically diverge, both for differential
equations \cite[Ch.~7]{Olver1974} and difference equations
\cite{Wong:1992:AESa}. This second solution underpins the connection formulas, the
Stokes phenomenon, and the comparison of the three resummations below.

\emph{(ii) Three resummations and their sectors of validity.} Each of the three
$q$-Laplace transforms in \S\ref{S:prelim} produces a solution
$y_e(z)=(\LqK{e}B_q)(z)$, $e=E,\theta,\lambda$, of \eqref{newexpandedequation}, all
sharing the asymptotic expansion \eqref{originaly} but valid on different sectors:
\begin{equation}\label{sectors}
    \begin{split}
        y_E(z)&\sim \sum_{n=0}^\infty c_n q^{-\binom{n}{2}} \left(-qz\right)^n,\qquad
        {\rm for~all~}\arg z,\\
        y_\theta(z)&\sim\sum_{n=0}^\infty c_n q^{-\binom{n}{2}} \left(-qz\right)^n,\qquad
        |\arg z|<2\pi,\\
        y_\lambda(z)&\sim \sum_{n=0}^\infty c_n q^{-\binom{n}{2}} \left(-qz\right)^n,\qquad
        |\arg (z/\lambda)|<\pi,
    \end{split}
\end{equation}
as $|z|\to 0$. Here $y_\lambda$ is single-valued, whereas $y_E$ and $y_\theta$ are
multi-valued.

\emph{(iii) Connection formulas and the Stokes multiplier.} The multi-valued
resummations satisfy
\begin{equation}\label{connectionsimple}
y_e(z)-y_e\left(z\expe^{2\pi\iunit}\right)=2\pi\iunit K(q)\,\secondsol\,u(z),\qquad e=E,\theta,
\end{equation}
so that crossing the negative real axis switches on a multiple of the
$q$-exponentially small second solution, with Stokes multiplier $K(q)$ determined in
\S\ref{S:vn}. The sectors of validity in \eqref{sectors} are a direct consequence of
\eqref{connectionsimple}. In the case of $y_E(z)$, the Stokes phenomenon switches on terms that 
will be $q$-exponentially small on all Riemann sheets. For $y_\theta(z)$, the switched-on
terms will be $q$-exponentially small until we encounter the poles of $1/\Thetaq{-z}$.

\emph{(iv) Exact resummations and error bounds.} In \S\ref{S:exact}, we obtain
convergent expansions in terms of so-called $q$-hyperterminant functions that represent the $y_e$
exactly. We then derive from them sharp error bounds for the truncated expansions of
\eqref{sectors}: simplest for $y_E$, and blowing up for $y_\lambda$ as
$|\arg(z/\lambda)|\to\pi$, where that resummation has poles by construction.

\emph{(v) Differences of the resummations.} In \S\ref{S:difference}, we prove
\begin{equation}\label{fulldifference}
\begin{split}
y_\theta(z)-y_E(z)&=K(q) \frac{u(z)}{\Thetaq{-q^2z}}P_q^{(c)}(z),\\
y_\theta(z)-y_\lambda(z)&=K(q) \frac{u(z)}{\Thetaq{-q^2z}}P_q^{(d)}(z;\lambda).
\end{split}
\end{equation}
That is, the Stokes multiplier times the $q$-exponentially small second solution times
\begin{align}
    P_q^{(c)}(z)    &=\frac{-2\pi}{\ln q}\sum_{n=1}^\infty \left(-1\right)^n\qhat^{n^2}
    \frac{\sin\left(2\pi n\frac{\ln z}{\ln q}\right)}{\sinh\left(n\ln \qhat\right)}\label{Pqc1},\\
    P_q^{(d)}(z;\lambda)&=\frac{\lambda\dThetaq{\lambda}}{\Thetaq{\lambda}}
    -\frac{\lambda\dThetaq{\lambda/z}}{z\Thetaq{\lambda/z}}+\frac{\ln z}{\ln q}.\label{Pqd1a}
\end{align}
These are \emph{universal} $q$-periodic functions, independent of $a(z), b(z)$.
The Fourier series \eqref{Pqc1} and \eqref{PqFourier}
explain why these differences are \emph{even smaller} than $q$-exponential size, as
illustrated numerically in \S\ref{S:Numerics}. For alternative representations of 
$P_q^{(c)}(z)$ and $P_q^{(d)}(z;\lambda)$, see \eqref{Pqd1} and \eqref{Pqc0}.

\emph{(vi) Asymptotics and applications.} In \S\ref{S:CH}, the sectors \eqref{sectors}
and the connection formula \eqref{connectionsimple} yield Cauchy--Heine
representations, late-coefficient asymptotics and exponentially-improved expansions.
These are illustrated for the basic hypergeometric series $\qphi{2}{0}{a,b}{-}{q}{z}$
(\S\ref{S:2phi0}), with the bilateral $\qpsi{0}{2}{-}{a,0}{q}{z}$ as a special case,
and in a non-rational numerical example (\S\ref{S:Numerics}), where a level-$1$
hyperasymptotic approximation can give us arbitrary accuracy!

\subsection{Related work and organisation}
In \cite{DiVizio2009}, the difference $y_\theta-y_\lambda$ is studied for the
$q$-Euler series; their result is equivalent to \eqref{v1diff}, and our presentation
of $P_q^{(d)}(z;\lambda)$ appears simpler. The difference of two
$\lambda$-resummations,
\begin{equation}\label{difflambdas}
\begin{split}
    P_q^{(d)}(z;\lambda_1)-P_q^{(d)}(z;\lambda_2)&=
    \frac{\lambda_1\dThetaq{\lambda_1}}{\Thetaq{\lambda_1}}
    -\frac{\lambda_1\dThetaq{\frac{\lambda_1}{z}}}{z\Thetaq{\frac{\lambda_1}{z}}}
    -\frac{\lambda_2\dThetaq{\lambda_2}}{\Thetaq{\lambda_2}}
    +\frac{\lambda_2\dThetaq{\frac{\lambda_2}{z}}}{z\Thetaq{\frac{\lambda_2}{z}}}\\
    &=\frac{\qpr{q}{q}{\infty}^3\Thetaq{-z}\Thetaq{\frac{-\lambda_1}{\lambda_2}}
    \Thetaq{\frac{-\lambda_1\lambda_2}{z}}}{\Thetaq{\lambda_1}\Thetaq{\lambda_2}
    \Thetaq{\frac{\lambda_1}{z}}\Thetaq{\frac{z}{\lambda_2}}},
\end{split}\end{equation}
recovers the quotient found in \cite[\S7.1]{Ramis2013}; the second identity is
Lemma~\ref{lemmaparfrac}.

The paper is organised as follows: \S\ref{S:general} reduces a general equation to
the canonical form \eqref{newexpandedequation}; \S\ref{S:secondsolution} constructs
the second solution; \S\ref{S:Wronskian}--\S\ref{S:vn} develop the Wronskian
relation, the $q$-hyperterminant functions $v_n$, and the Stokes multiplier;
\S\ref{S:exact} gives exact resummations and error bounds; \S\ref{S:difference} determines
the differences of the resummations; \S\ref{S:CH} treats Cauchy--Heine asymptotics;
\S\ref{S:2phi0} examines the hypergeometric example; and
\S\ref{S:Numerics} considers the numerics. Appendices~\ref{Apptheta}--\ref{AppLemma} collect
theta function identities and the coefficient bounds, respectively.

\section{The \texorpdfstring{$q$}{q}-Laplace transforms}\label{S:prelim}

\subsection{Two weight functions, one moment problem}
To invert \eqref{qBorel}, we need a weight $\K{t}$ with moments
$\int_0^\infty t^n\K{t}\,\intd t=q^{-\binom{n+1}{2}}$. This is the indeterminate
Stieltjes moment problem of the Stieltjes--Wigert polynomials, so the weight is not
unique. Following \cite{Christiansen2003}, the two canonical choices are given in \eqref{X}.

Although weights $\Cq\Eq{t}$ and $\frac{-1}{\ln(q)\Thetaq{t}}$ share these moments, they
are different functions, related by
\begin{equation}\label{defPq}
    \Cq\Eq{t}=\frac{-1}{\ln(q)\Thetaq{t}}\Pq{t},
\end{equation}
with $\Pq{t}$ $q$-periodic: $\Pq{qt}=\Pq{t}$. Appendix~\ref{Apptheta} shows via
\eqref{PqFourier} that $\Pq{t}-1$ is extremely small for $q$ bounded away from $0$.

\subsection{The three \texorpdfstring{$q$}{q}-Laplace transforms}
Many $q$-Laplace transforms appear in the literature, e.g.\ \cite{Tahara2017}; we
require the property $t^n\mapsto q^{-\binom{n}{2}}z^n$. The two continuous
transforms \cite{Ramis:1993:SDA,Zhang1999,Zhang2000} are
\begin{equation}\label{qLaplace1}
    (\LqK{E} B)(z)=\Cq\int_0^\infty B(t)\Eq{\frac{t}{z}}\frac{\intd t}{t}, \qquad (\LqK{\theta} B)(z)=\frac{-1}{\ln q}\int_0^\infty \frac{B(t)}{\Thetaq{\frac{t}{z}}}\frac{\intd t}{t},
\end{equation}
and the discrete transform \cite{Zhang2002} is
\begin{equation}\label{qLaplace2}
    (\LqK{\lambda} B)(z)=\sum_{n=-\infty}^{\infty}\frac{B(q^n\lambda)}{\Thetaq{\frac{q^n\lambda}{z}}}
    =\frac{1}{1-q}\int_0^\infty \frac{B(t\lambda)}{\Thetaq{\frac{t\lambda}{z}}}\frac{\intd_q t}{t},
\end{equation}
for $z\in\mathbb{C}\backslash\{-q^n\lambda;n\in\mathbb{Z}\}$, using the $q$-integral
notation \cite[\href{https://dlmf.nist.gov/17.2.E48}{17.2.48}]{NIST:DLMF}. There is
no discrete analogue built from $\Eq{t}$. Since $B(t)$ has singularities in the
$t$-plane, the nodes $q^n\lambda$ must avoid them; when $B$ has poles at $t=-q^{-m}$,
we require $\lambda\in\mathbb{C}\backslash\{-q^n;n\in\mathbb{Z}\}$.

Finally, the continuous transforms arise from the discrete one by integrating over
the free parameter $\lambda$:
\begin{equation}\label{Intlamda}
    (\LqK{\theta} B)(z)=\frac{-1}{\ln q}\int_q^1 (\LqK{\lambda} B)(z)\,\frac{\intd \lambda}{\lambda},
\end{equation}
\cite[Theorem 4.14]{DiVizio2009}, and
\begin{equation}\label{IntlamdaE}
    (\LqK{E} B)(z)=\frac{-1}{\ln q}\int_q^1 (\LqK{\lambda} B)(z)\,
    \Pq{\frac{\lambda}{z}}\frac{\intd \lambda}{\lambda}.
\end{equation}

\subsection{The \texorpdfstring{$q$}{q}-Laplace transform table}
In Table \ref{tab:qLaplace}, we display the $q$-Laplace transform pairs: given $f(t)$ 
(Borel side, left), the table lists $F(z)=(\LqK{\ee} f)(z)$ (Laplace side, right). Note that
property (3) is slightly more subtle for the discrete $q$-Laplace transform:
$(\LqK{\lambda} f(q/t))(z)=(\LqK{1/\lambda} f(t))(1/z)$.
\begin{table}[ht]
\begin{center}
\begin{tabular}{|c|c|c|}
\hline
&$f(t)$ & $F(z)=(\LqK{\ee} f)(z)$ \\
\hline

(1)&$t^nf(t)$ & $\displaystyle q^{-\binom{n}{2}} z^nF(zq^{-n})$\\

(2)&$f(qt)$ & $\displaystyle F(qz)$\\

(3)&$f(q/t)$ & $\displaystyle F(1/z)$\\

\hline

(4)&$t^n$ & $\displaystyle q^{-\binom{n}{2}} z^n$\\

(5)&$\qpr{-at}{q}{\infty}$ & $\displaystyle\frac1{\qpr{az}{q}{\infty}}$\\

(6)&$\displaystyle \qpr{-bq/t}{q}{\infty}$ & $\displaystyle \frac{1}{\qpr{b/z}{q}{\infty}}$\\

(7)&$\displaystyle \qpr{-at}{q}{\infty}\qpr{-bq/t}{q}{\infty}$ & 
$\displaystyle \frac{\qpr{ab}{q}{\infty}}{\qpr{az}{q}{\infty}\qpr{b/z}{q}{\infty}}$\\

(8)&$\displaystyle \frac{\qpr{-at}{q}{\infty}}{\qpr{-bt}{q}{\infty}}$ & 
$\displaystyle \UqK{\ee}{0}{\frac{a}{b}}{bz}$,~~$\ee=E,\theta$\\

(9)& $\displaystyle \qphi{1}{1}{q^{-n}}{0}{q}{-tq^n}$ & $\qpr{z}{q}{n}$ \\

(10)& $\displaystyle \qphi{r}{s+1}{a_1,\ldots,a_r}{b_1,\ldots,b_s,0}{q}{-t}$ & $\displaystyle \qphi{r}{s}{a_1,\ldots,a_r}{b_1,\ldots,b_s}{q}{z}$ \\

\hline
\end{tabular}
\caption{Rules (1)--(3) are operational
identities; (4)--(10) are explicit pairs. For (5)--(7) we require $b<z<a^{-1}$, and
(7) follows from \cite[Thm 6.3]{Cohl2025}. The function $\UqK{\ee}{0}{\frac{a}{b}}{bz}$
in (8) is defined in \eqref{2phi0BorelLaplaceK}.}

\label{tab:qLaplace}
\end{center}
\end{table}%


\section{General second order linear $q$-difference equations}\label{S:general}
A general second order linear $q$-difference equation, for which formal solutions
are divergent, is of the form
\begin{equation}\label{originaleq}
xw(xq^{-2})+f(x)w(xq^{-1})-g(x)w(x)=0,
\end{equation}
where $f(x)=\sum\limits_{m=0}^\infty f_m x^m$, $g(x)=\sum\limits_{m=0}^\infty g_m x^m$ have a 
non-zero radius of convergence. The $q$-difference equation for the divergent basic hypergeometric
series ${}_{2}\phi_{0}$ is of this form with linear functions $f(x)=q-(a+b)x$ and $g(x)=q-abx$.

We consider the lowest order case of $f(x)$ and $g(x)$ where $f_0\neq0$, $g_0\neq0$. First we look for a balance
between the second and third term in \eqref{originaleq}. In the first instance,
we have to solve $f_0 \widetilde{w}(xq^{-1})-g_0 \widetilde{w}(x)=0$ with solution 
$\widetilde{w}(x)=\left(f_0/g_0\right)^{\frac{\ln x}{\ln q}}p(x)$, 
in which $p(x)$ is $q$-periodic, that is, $p(qx)=p(x)$ for all $x$.
We make the substitution,
\begin{equation}
    w(x)= \left(\sfrac{f_0}{g_0}\right)^{\frac{\ln x}{\ln q}}v(x),
\end{equation}
and obtain the equation,
\begin{equation}
    \frac{g_0}{f_0^2}xv(xq^{-2})+\frac{f(x)}{f_0} v(xq^{-1})-\frac{g(x)}{g_0} v(x)=0.
\end{equation}
It then makes sense to introduce a new variable $z=\frac{g_0}{f_0^2}x$ and let $y(z)=v(x)$. Then \eqref{newexpandedequation} with \eqref{abz} hold.


\section{The second solution}\label{S:secondsolution}

To find the second solution, we look for a balance between the first and second term in 
\eqref{newexpandedequation}. In this instance, we first have to solve 
$z\tilde{y}(zq^{-2})+\tilde{y}(zq^{-1})=0$ or $\tilde{y}(qz)=-q^2z\tilde{y}(z)$. 
Noting that $\K{qz}=z\K{z}$, we find $\tilde{y}(z)=\K{\expe^{\pi\iunit}q^2z}$. 
Letting $y(z)=\K{\expe^{\pi\iunit}q^2z}u(z)$, we have,
\begin{equation}\label{uequation}
    u(zq^{-2})-a(z)u(zq^{-1})-qzb(z)u(z)=0.
\end{equation}
This equation then has a formal solution of the form,
\begin{equation} \label{uexpansion}
    u(z)=\sum_{n=0}^\infty d_n z^n,
\end{equation}
with
\begin{equation}
    d_n=\frac{q^{2n}}{1-q^n}\sum_{k=0}^{n-1}q^{-k}d_k\left(a_{n-k}+q^{k+1}b_{n-k-1}\right).
\end{equation}
We take $d_0=1$.
The expansion \eqref{uexpansion} has a non-zero radius of convergence. 
This is proved in Appendix \ref{AppLemma}.


\section{Wronskian Relation}\label{S:Wronskian}

Let $y(z)$ be a solution of \eqref{newexpandedequation} with the property 
$y(z)\sim 1$ as $z\to 0$, and let $y_u(z)=\K{\expe^{\pi\iunit}q^2z}u(z)$ be 
the second solution defined in the previous section. We define the Wronskian via
\begin{equation}
    W(z)=y(z)y_u\left(zq^{-1}\right)-y\left(zq^{-1}\right)y_u(z).
\end{equation}
When we consider $z\to 0$, we obtain that $W(z)\sim \K{\expe^{\pi\iunit}qz}$.
The Wronskian satisfies
\begin{equation}
    zW\left(zq^{-1}\right)=-b(z)W(z).
\end{equation}
Then
\begin{equation}
    \frac{W(z)}{\K{\expe^{\pi\iunit}qz}}=\frac{W(z q^n)}{\K{\expe^{\pi\iunit}qz q^n}}
    \prod_{j=1}^n b(zq^j),
\end{equation}
from which we obtain 
\begin{equation}
    \frac{W(z)}{\K{\expe^{\pi\iunit}qz}}=\prod_{j=1}^\infty b(zq^j)=\beta(z).
\end{equation}
That is,
\begin{equation} \label{wronskianu}
    y(z)u\left(zq^{-1}\right)+zqy\left(zq^{-1}\right)u(z)
    =\beta(z).
\end{equation}
Recall that we assume that $b(z)$ is analytic near $z=0$ and $b(0)=1$. 
It follows that $\beta(z)$ is also
analytic near $z=0$.
To obtain a recurrence relation for the terms of $\beta(z)$, let
\begin{equation}
    \beta(z)
    =\sum_{n=0}^\infty \beta_n z^n,
\end{equation}
with $\beta_0=1$. Using $\beta\left(zq^{-1}\right)=b(z)\beta(z)$, we then have, 
for $n\geq1$,
\begin{equation}
    \beta_n=\frac{q^n}{1-q^n}\sum_{k=0}^{n-1} \beta_k b_{n-k}.
\end{equation}

We may write the Wronskian relation \eqref{wronskianu} as
\begin{equation}\label{talpha}
    \frac{y(qz)}{u(qz)}+zq^2\frac{y(z)}{u(z)}=\frac{\beta(qz)}{u(qz)u(z)}
    =\sum_{n=0}^\infty \alpha_n z^n,
\end{equation}
and solve
\begin{equation} \label{vnequation}
    v_n(qz)+zq^2v_n(z)=z^n,
\end{equation}
in terms of the $q$-hyperterminant function.
In the next section, we will obtain several solutions of \eqref{vnequation} and show that
\begin{equation}\label{vnexpansion}
    y(z)=u(z)\sum_{n=0}^\infty\alpha_nv_n(z)
\end{equation}
is convergent. A formal solution of \eqref{vnequation} is
\begin{equation}\label{vnformal}
    v_n(z)\sim \left(-1\right)^{n}q^{n(n-5)/2}\sum_{m=n}^\infty q^{-\binom{m}{2}} \left(-qz\right)^{m}.
\end{equation}
Reversing the derivation shows that \eqref{vnexpansion} indeed solves \eqref{newexpandedequation}.

We have $\alpha_0=1$, and obtain from \eqref{talpha}, the recurrence relation
\begin{equation}
    \alpha_n={\beta}_nq^n-\sum_{m=0}^{n-1}\alpha_m\sum_{k=0}^{n-m}{d}_k{d}_{n-m-k}q^k,
\end{equation}
for $n\geq1$.

\section{Solving the $v_n$ equation}\label{S:vn}

We have presented \eqref{vnformal} in such a way that the infinite series is obviously the
$q$-Laplace transform of $\frac{t^n}{1+t}$. Hence, three solutions of \eqref{vnequation} are
\begin{align}
v_{n,e}(z)&=q^{\frac{1}{2}n(n-5)}\int_0^\infty \frac{t^{n-1}\K{\frac{t}{qz}}}{1+t}\intd t,
    \qquad e=E,\theta,\label{vne}\\
    v_{n,\lambda}(z)&=q^{\frac{1}{2}n(n-5)}\sum_{m=-\infty}^\infty
    \frac{\lambda^n q^{mn}}{\left(1+\lambda q^m\right)\Thetaq{\frac{\lambda q^m}{qz}}}.\label{vnl}
\end{align}
When we combine these representations with \eqref{Fprop}, it is easy to show that
\begin{equation} \label{vnmqrecurrence}
    z^m v_{n,e}(z)=q^{-m^2-(n-1)m}v_{n+m,e}(zq^m),\qquad e=E,\theta,\lambda.
\end{equation}

We now continue with the cases $e=E,\theta$.
Our expansion \eqref{vnexpansion} for $y(z)$, with $u(z)=\sum\limits_{m=0}^\infty d_m z^m$, 
is then given by
\begin{equation} \label{yvnexpression}
\begin{split}
    y_e(z)&=\sum_{m=0}^\infty d_m z^m \sum_{n=0}^\infty \alpha_n v_{n,e}(z)
    =\sum_{m=0}^\infty \sum_{n=0}^\infty d_m \alpha_n q^{-m^2-(n-1)m}
    v_{n+m,e}(zq^m)\\
    &=\sum_{m=0}^\infty \sum_{n=0}^\infty d_m \alpha_n 
    q^{-m^2-(n-1)m+\frac{1}{2}(n+m)(n+m-5)}\int_0^\infty \frac{\tau^{n+m-1}
    \K{\frac{\tau}{q^{m+1}z}}}{1+\tau}\intd\tau\\
    &=\sum_{m=0}^\infty \sum_{n=0}^\infty d_m \alpha_n 
    q^{\frac{1}{2}(n+m)^2-\frac{3}{2}(n+m)-1}\int_0^\infty \frac{t^{n+m-1}
    \K{\frac{t}{z}}}{t+q^{-m-1}}\intd t\\
    &=\sum_{m=0}^\infty \sum_{\ell=m}^\infty d_m \alpha_{\ell-m} 
    q^{\frac{1}{2}\ell^2-\frac{3}{2}\ell-1}\int_0^\infty \frac{t^{\ell-1}
    \K{\frac{t}{z}}}{t+q^{-m-1}}\intd t.
\end{split}
\end{equation}
The location and degree of the poles in this expression for $y(z)$ agree with those found from the numerical calculation of the Pad\'e approximants of the right-hand side of \eqref{qBorel}.

Hence, we have for $y_e(z)$, the $q$-Laplace transform representation
\begin{equation}\label{BorelLaplace}
    y_e(z)=\int_0^\infty B_q(t) \K{\frac{t}{z}}\frac{dt}{t},
\end{equation}
in which the  $q$-Borel transform $B_q(t)$ has representation
\begin{equation} \label{qborel}
\begin{split}
    B_q(t)&=\sum_{m=0}^\infty \sum_{\ell=m}^\infty d_m \alpha_{\ell-m}q^{\frac{1}{2}\ell^2-\frac{3}{2}\ell-1}\frac{t^\ell}{t+q^{-m-1}}\\
    &=\sum_{\ell=0}^\infty t^\ell q^{\frac{1}{2}\ell^2-\frac{3}{2}\ell-1}\sum_{m=0}^\ell 
    \frac{d_m \alpha_{\ell-m}}{t+q^{-m-1}}.
\end{split}
\end{equation}

We take $\left|t+q^{-m-1}\right|\geq\delta$ for all $m$.
From earlier in the paper, we obtain that both $\sum\limits_{n=0}^\infty \left|d_n\right| z^n$
and $\sum\limits_{n=0}^\infty \left|\alpha_n\right| z^n$ are analytic near $z=0$. Hence, the
product of these two functions is also analytic at $z=0$, and we can deduce that
there exist positive constants $M$ and $D$ such that
$\sum\limits_{m=0}^\ell\left|d_m\alpha_{\ell-m}\right|\leq M D^\ell\leq
M\frac{D^\ell}{\qpr{q}{q}{\ell}}$ for all $\ell$. We estimate
\begin{equation} \label{qborel2}
    \left|B_q(t)\right|\leq\frac{M}{q\delta}\sum_{\ell=0}^\infty 
    \frac{\left(|t|D\right)^\ell q^{\frac{1}{2}\ell^2-\frac{3}{2}\ell}}{\qpr{q}{q}{\ell}}
    =\frac{M}{q\delta}\qpr{-D|t|/q}{q}{\infty}.
\end{equation}
When we use \cite[(2.3)]{Joshi2025}, we obtain \eqref{Bqest}.
Hence, for $|z|$ small enough and the path of integration bounded away from the poles, the
integral in \eqref{BorelLaplace} converges.

Expansion \eqref{qborel} has information about the residues of $B_q(t)$. The normal $q$-Borel
transform of \eqref{originaly} is 
\begin{equation} \label{normalqBorel}
    B_q(t)=\sum_{n=0}^\infty c_n \left(-qt\right)^n,\qquad |tq|<1.
\end{equation}
We can use the Cauchy integral representation for the coefficients in this expansion, but replace
$B_q(t)$ by \eqref{qborel}. After several steps, we will obtain the relation
\begin{equation}\label{cndn}
    c_n=\sum_{m=0}^n d_m q^{-\binom{m}{2}+(n-2)m}
    \sum_{\ell=0}^{n-m}\alpha_{\ell}q^{\frac{1}{2}\ell(\ell-5)}\left(-1\right)^{\ell+m}.
\end{equation}

We obtain from integral representation \eqref{vne}, the analytic connection formulas
\begin{equation}
    v_{n,e}(z\expe^{2\pi \iunit})-v_{n,e}(z)
    =2\pi\iunit \left(-1\right)^{n-1}q^{\frac{1}{2}n(n-5)}\K{\expe^{\pi \iunit}q^2z}.
\end{equation}
Hence, the analytic connection formulas for our solutions are
\begin{equation}\label{connection}
\begin{split}
    y_e(z\expe^{2\pi \iunit})-y_e(z)&=
    u(z)\sum_{n=0}^\infty \alpha_n\left(v_{n,e}(z\expe^{2\pi \iunit})-v_{n,e}(z)\right)\\
   &=-2\pi\iunit K(q)\K{\expe^{\pi \iunit}q^2z}u(z),
\end{split}
\end{equation}
with Stokes multiplier
\begin{equation}\label{Kq}
    K(q)=\sum_{n=0}^\infty \alpha_n\left(-1\right)^{n}q^{\frac{1}{2}n(n-5)}.
\end{equation}
Thus, when $z$ crosses the negative real axis, the Stokes phenomenon switches on 
a multiple of the second solution of our $q$-difference equation.


\section{Exact resummation and error-bounds}\label{S:exact}
The final representation in \eqref{yvnexpression} can be presented as
\begin{equation}\label{exactresummation}
    y_e(z)=\sum_{\ell=0}^\infty\sum_{m=0}^\ell d_m \alpha_{\ell-m}z^m
    v_{\ell-m,e}(z),\qquad e=E,\theta,\lambda.
\end{equation}
The case $e=\lambda$ is very similar to the other two cases. This double series converges
absolutely for $|z|$ small enough, and it can be seen as an exact resummation as long as we can evaluate the $q$-hyperterminant functions $v_{n,e}(z)$. It is easy to check that we have
the recurrence relation
\begin{equation}\label{vnrecurrence}
    v_{n+1,e}(z)=q^{-2}z^n-q^{n-2}v_{n,e}(z).
\end{equation}
Hence, if we accept $v_{0,e}(z)$ as an acceptable approximant, then we have the exact resummation
\eqref{exactresummation}.

When we write \eqref{vnrecurrence} as $v_{n,e}(z)=q^{-n}z^n-q^{2-n}v_{n+1,e}(z)$, we obtain 
for the finite sum in \eqref{exactresummation} that 
$\sum\limits_{m=0}^\ell d_m \alpha_{\ell-m}z^mv_{\ell-m,e}(z)
=\bigO\left(z^\ell\right)$ as $z\to 0$. We combine these observations and obtain for the remainder
\begin{equation}\label{Remainder}
    R_{N,e}(z)=y_e(z)-\sum_{n=0}^{N-1}c_n q^{-\binom{n}{2}} 
    \left(-qz\right)^n,
\end{equation}
the convergent expansion
\begin{equation}\label{RNresummation}
\begin{split}
    R_{N,e}(z)&=\sum_{\ell=0}^{N-1}\sum_{m=0}^\ell \left(-1\right)^{\ell-N}
    q^{-\frac12(N-\ell)(N+\ell-2m-5)}d_m \alpha_{\ell-m}z^m
    v_{N-m,e}(z)\\
    &\qquad+\sum_{\ell=N}^\infty\sum_{m=0}^\ell d_m \alpha_{\ell-m}z^m
    v_{\ell-m,e}(z).
\end{split}\end{equation}

\begin{theorem} Let $M$ and $D$ be positive constants such that
$\sum\limits_{m=0}^\ell \left|d_m 
\alpha_{\ell-m}\right|\leq M D^\ell\leq \frac{M D^\ell}{\qpr{q}{q}{\ell}}$
for all $\ell$. We take $|z|<q/D$ and $\lambda>0$. Then we have the error-bound
\begin{equation}\label{RNbound0}
    \left|R_{N,e}(z)\right|\leq 
    C M\left|z\right|^N\left(q^{-\frac12 N(N-3)}\qpr{-Dq^{-2}}{q}{\infty}
    +\frac{\left(D/q\right)^N}{1-\frac{|z|D}{q}}\right),
\end{equation}
with $C=1$ in the half-plane $|\arg z|\leq\frac12\pi$ for $e=E,\theta,\lambda$, and
in the sectors $\xi=|\arg z|\in(\frac12\pi,\pi)$, we have
\begin{equation}\label{RNboundC}
C=\begin{cases}
\exp\left(\frac{-\zeta^2}{2\ln q}\right),& e=E,\\
\frac{\exp\left({\frac{-\zeta^2}{2\ln q}}\right)}{\qpr{\qhat^2}{\qhat^2}{\infty}^3}\left(\frac{1+\sqrt{\qhat}}{1-\sqrt{\qhat}}\right),& e=\theta,\\
\frac1{\left|\sin\xi\right|},& e=\lambda,
\end{cases}
\end{equation}
in which $\zeta=|\arg z|-\frac12\pi$.
\end{theorem}
The first term on the right-hand side of \eqref{RNbound0} dominates for large $N$. In the central
sector, $|\arg z|\leq\frac12\pi$, all three resummations have the same error-bound, and in the
sectors $|\arg z|\in(\frac12\pi,\pi)$, the error-bound for the case $e=E$ is the smallest, the one for
the case $e=\theta$ is slightly larger, and the one for the case $e=\lambda$ blows up when
we approach the boundary $|\arg z|=\pi$, which is where the resummation has poles by construction.

\begin{proof}
In the proof, we will show that $\left|v_{n,e}(z)\right|\leq C q^{-n}\left|z\right|^n$.
Once we have this result, we can use it in \eqref{RNresummation} and obtain
\begin{equation}\label{RNbound1}
\begin{split}
    \left|R_{N,e}(z)\right|&\leq C\left|z\right|^N\sum_{\ell=0}^{N-1}
    q^{-\frac12 N(N-3)+\frac12\ell(\ell-5)}\sum_{m=0}^\ell 
    \left|d_m \alpha_{\ell-m}\right|\\
    &\qquad+C\sum_{\ell=N}^\infty\left|z\right|^\ell q^{-\ell}\sum_{m=0}^\ell 
    \left|d_m \alpha_{\ell-m}\right|\\
    &\leq CM\left|z\right|^Nq^{-\frac12 N(N-3)}\sum_{\ell=0}^{\infty}
    \frac{q^{\frac12\ell(\ell-5)}D^\ell}{\qpr{q}{q}{\ell}}
    +CM\sum_{\ell=N}^\infty\left|z\right|^\ell q^{-\ell}D^\ell,
\end{split}\end{equation}
and \eqref{RNbound0} follows.

We start with the cases $e=E,\theta$ and $\Re z\geq 0$. We do need the bound
\begin{equation}\begin{split}
    \left|v_{n,e}(z)\right|&=q^{\frac{1}{2}n(n-3)}\left|z\right|^n\left|\int_0^\infty 
    \frac{t^{n-1}\K{t}}{1+tqz}\intd t\right|\\
    &\leq q^{\frac{1}{2}n(n-3)}\left|z\right|^n\int_0^\infty 
    t^{n-1}\K{t}\intd t=q^{-n}\left|z\right|^n,
\end{split}\end{equation}
that is, $C=1$. Now we take $z=|z|\expe^{\iunit(\zeta+\frac\pi2)}$ with
$\zeta\in(0,\frac12\pi)$. In the case $e=E$, we have
\begin{equation}\begin{split}
    v_{n,E}(z)&=\Cq q^{\frac{1}{2}n(n-3)}z^n\int_0^{\infty \expe^{-\iunit\zeta}}
    \frac{\tau^{n-1}\Eq{\tau}}{1+\tau qz}\intd \tau\\
    &= \Cq \expe^{-\iunit n\zeta}q^{\frac{1}{2}n(n-3)}z^n\int_0^{\infty}
    \frac{t^{n-1}\Eq{t\expe^{-\iunit\zeta}}}{1+t q|z|\iunit}\intd t.
\end{split}\end{equation}
Hence,
\begin{equation}\begin{split}
    \left|v_{n,E}(z)\right|&=\Cq q^{\frac{1}{2}n(n-3)}\left|z\right|^n
    \exp\left(\frac{-\zeta^2}{2\ln q}\right)\int_0^{\infty \expe^{-\iunit\zeta}}
    \frac{t^{n-1}\Eq{t}}{\left|1+\iunit t q|z|\right|}\intd t\\
    &\leq q^{-n}\left|z\right|^n\exp\left(\frac{-\zeta^2}{2\ln q}\right),
\end{split}\end{equation}
and the first case of \eqref{RNboundC} follows. We use \eqref{defPq} and obtain in a similar manner
\begin{equation}\label{vnthetabound}
    v_{n,\theta}(z)=\frac{-\expe^{-\iunit n\zeta}}{\ln q}q^{\frac{1}{2}n(n-3)}z^n\int_0^{\infty}
    \frac{t^{n-1}\Eq{t\expe^{-\iunit\zeta}}}{(1+\iunit t q|z|)\Pq{t\expe^{-\iunit\zeta}}}\intd t.
\end{equation}
We use \eqref{PqReciprocal} to obtain a bound for $1/|\Pq{t\expe^{-\iunit\zeta}}|$. We have
\begin{equation}
    \left|\cos\left(\frac{2\pi n}{\ln q}(-\iunit\zeta+\ln t)\right)\right|
    \leq\expe^{-\frac{2\pi n|\zeta|}{\ln q}}
    \leq\expe^{-\frac{\pi^2n}{\ln q}}=\qhat^{-\frac{n}{2}}\quad\text{for}\quad\zeta\in\left(-\tfrac{\pi}{2},\tfrac{\pi}{2}\right),
\end{equation}
and note that $\widetilde{a}_n\leq1$ since the terms in expansion \eqref{tan} are alternating and monotonically decreasing in magnitude. Using these, we have
\begin{equation}
\begin{split}
    \frac{1}{|\Pq{t\expe^{-\iunit\zeta}}|}&\leq\frac{1}{\qpr{\qhat^2}{\qhat^2}{\infty}^3}\left(|\widetilde{a}_0|+2\sum_{n=1}^\infty|\widetilde{a}_n|\qhat^n
    \left|\cos\left(\frac{2\pi n}{\ln q}(-\iunit\zeta+\ln t)\right)\right|\right)\\
    &\leq\frac{1}{\qpr{\qhat^2}{\qhat^2}{\infty}^3}\left(1+2\sum_{n=1}^\infty\qhat^{\frac{n}{2}}\right)
    =\frac{1}{\qpr{\qhat^2}{\qhat^2}{\infty}^3}\frac{1+\sqrt{\qhat}}{1-\sqrt{\qhat}}.
\end{split}
\end{equation}
We combine this with \eqref{vnthetabound}, take absolute values and the second case 
of \eqref{RNboundC} follows.

For the discrete $q$-Laplace transform, we cannot use these substitutions and rotations of the
contour of integration. We combine \eqref{PartFrac1a} with the observation that
$\sum\limits_{m\in\mathbb{Z}}\frac1{\Thetaq{\mu q^m}}=1$ and obtain
\begin{equation}\begin{split}
    &q^n z^{-n}v_{n,\lambda}(z)=z^{-n}q^{\frac{1}{2}n(n-3)}\sum_{m=-\infty}^\infty
    \frac{\lambda^n q^{mn}}{\left(1+\lambda q^m\right)\Thetaq{\frac{\lambda q^m}{qz}}}\\
    &\qquad=\sum_{m=-\infty}^\infty
    \frac{1}{\left(1+\lambda q^m\right)\Thetaq{\frac{\lambda q^{n+m}}{qz}}}
    =1-\lambda\sum_{m=-\infty}^\infty
    \frac{1}{\left(\lambda+ q^{-m}\right)\Thetaq{\frac{\lambda q^{n+m}}{qz}}}\\
   &\qquad=1-\lambda\sum_{m=-\infty}^\infty
    \frac{1}{\left(\frac{\lambda q^{n}}{qz}+ q^{-m}\right)\Thetaq{\lambda q^m}}
    =\sum_{m=-\infty}^\infty
    \frac{1}{\left(1+\frac{qz}{\lambda q^{n+m}}\right)\Thetaq{\lambda q^{m+1}}}.
\end{split}\end{equation}
In the case that $\Re z\geq 0$, we use that $\left|1+\frac{qz}{\lambda q^{n+m}}\right|^{-1}\leq 1$
and obtain $C=1$.

Finally, when $\xi=\left|\arg z\right|\in(\frac\pi2,\pi)$, we observe that the function $r\to\left|\frac{1}{1+rz}\right|^2=\frac{1}{r^2|z|^2+1+2r|z|\cos\xi}$
has a maximum when $r|z|=-\cos\xi$. It follows that 
$\left|\frac{1}{1+rz}\right|\leq \frac{1}{|\sin\xi|}$.
\end{proof}


\section{The differences of the $q$-Laplace transforms}\label{S:difference}

We are going to determine $v_{n,\theta}(z)-v_{n,\lambda}(z)$, but we first make the observation that
we can focus on the case $n=1$ because it follows from definitions \eqref{vne} and \eqref{vnl} that
$v_{0,e}(z)+q^2v_{1,e}(z)=1$ and for $n\geq1$,
\begin{equation}
    \left(-1\right)^{n-1} q^{\frac12n(5-n)}v_{n,e}(z)=q^2v_{1,e}(z)+\sum_{m=1}^{n-1}
    \left(-qz\right)^m q^{-\binom{m}{2}}.
\end{equation}
Thus
\begin{equation}
    \left(-1\right)^{n-1} q^{\frac12n(5-n)}\left(v_{n,\theta}(z)-v_{n,\lambda}(z)\right)
    =q^2\left(v_{1,\theta}(z)-v_{1,\lambda}(z)\right).
\end{equation}
We use \eqref{PartFrac1} twice and partial fractions,
\begin{align}
    &q^2\left(v_{1,\theta}(z)-v_{1,\lambda}(z)\right)=\sum_{n=-\infty}^\infty
    \left(\frac{-1}{\ln q}\int_{q^{n+1}}^{q^n} \frac{\intd t}{\Thetaq{\frac{t}{qz}}(1+t)}-
    \frac{q^n\lambda}{\Thetaq{\frac{q^n\lambda}{qz}}(q^n\lambda+1)}\right)\\
    &\qquad =\sum_{n=-\infty}^\infty
    \left(\frac{-1}{\ln q}\int_{q}^{1} \frac{q^n\intd t}{\Thetaq{\frac{q^nt}{qz}}(q^nt+1)}-
    \frac{q^n\lambda}{\Thetaq{\frac{q^n\lambda}{qz}}(q^n\lambda+1)}\right)\notag\\
    &\qquad =\sum_{n=-\infty}^\infty
    \left(\int_{0}^{1} \frac{q^{n+\tau}\intd \tau}{\Thetaq{\frac{q^{n+\tau}}{qz}}(q^{n+\tau}+1)}-
    \frac{q^n\lambda}{\Thetaq{\frac{q^n\lambda}{qz}}(q^n\lambda+1)}\right)\notag\\
    &\qquad =\sum_{n=-\infty}^\infty
    \sum_{m=-\infty}^\infty 
    \frac{\left(-1\right)^mq^{\binom{m}{2}}}{\qpr{q}{q}{\infty}^3}
    \left(\int_{0}^{1} \frac{q^\tau\intd \tau}{\left(\frac{q^\tau}{qz}+q^{-m-n}\right)(q^{n+\tau}+1)}-
    \frac{\lambda}{\left(\frac{\lambda}{qz}+q^{-m-n}\right)(q^n\lambda+1)}\right)\notag\\
    &\qquad =\sum_{n=-\infty}^\infty
    \sum_{m=-\infty}^\infty 
    \frac{\left(-1\right)^mq^{\binom{m}{2}}}{\qpr{q}{q}{\infty}^3(q^{-m}-\frac1{qz})}
    \left(\int_{0}^{1} \frac{q^\tau\intd \tau}{q^\tau+q^{-n}}-
    \frac{\lambda}{\lambda+q^{-n}}\right)\notag\\
    &\qquad\qquad-\sum_{n=-\infty}^\infty
    \sum_{m=-\infty}^\infty 
    \frac{\left(-1\right)^mq^{\binom{m}{2}}}{\qpr{q}{q}{\infty}^3(q^{-m}-\frac1{qz})}
    \left(\int_{0}^{1} \frac{q^\tau\intd \tau}{q^\tau+zq^{1-m-n}}-
    \frac{\lambda}{\lambda+zq^{1-m-n}}\right)\notag\\
    &\qquad =\sum_{n=-\infty}^\infty
    \sum_{m=-\infty}^\infty 
    \frac{\left(-1\right)^mq^{\binom{m}{2}}}{\qpr{q}{q}{\infty}^3(q^{-m}-\frac1{qz})}
    \left(\int_{0}^{1} \frac{q^\tau\intd \tau}{q^\tau+q^{-n}}-
    \frac{\lambda}{\lambda+q^{-n}}\right)\notag\\
    &\qquad\qquad-\sum_{n=-\infty}^\infty
    \sum_{m=-\infty}^\infty 
    \frac{\left(-1\right)^mq^{\binom{m}{2}}}{\qpr{q}{q}{\infty}^3(q^{-m}-\frac1{qz})}
    \left(\int_{0}^{1} \frac{q^\tau\intd \tau}{q^\tau+zq^{1-n}}-
    \frac{\lambda}{\lambda+zq^{1-n}}\right)\notag\\
    &\qquad =\frac1{\Thetaq{\frac{-1}{qz}}}\sum_{n=-\infty}^\infty
    \left(\int_{0}^{1} \frac{q^\tau\intd \tau}{q^\tau+q^{-n}}-\int_{0}^{1} \frac{q^\tau\intd \tau}{q^\tau+zq^{1-n}}
    +\frac{\lambda}{\lambda+zq^{1-n}}
    -\frac{\lambda}{\lambda+q^{-n}}\right)\notag\\
    &\qquad =\frac1{\Thetaq{-q^2z}}\sum_{n=-\infty}^\infty
    \left(\frac1{\ln q}\ln\left(\frac{(q+q^{-n})(1+zq^{1-n})}{(1+q^{-n})(q+zq^{1-n})}\right)
   +\frac{\lambda}{\lambda+zq^{1-n}}
    -\frac{\lambda}{\lambda+q^{-n}}\right).\notag
\end{align}
The sum of the first term is a telescoping sum, and for the final two terms, we use 
\eqref{lemmaparfrac1}:
\begin{equation}\label{v1diff}
\begin{split}
    q^2\left(v_{1,\theta}(z)-v_{1,\lambda}(z)\right)&=\frac1{\Thetaq{-q^2z}}\left(-\frac{\ln (qz)}{\ln q}
    +\frac{\frac{\lambda}{qz}\dThetaq{\frac{\lambda}{qz}}}{\Thetaq{\frac{\lambda}{qz}}}
    -\frac{\lambda\dThetaq{\lambda}}{\Thetaq{\lambda}}
    \right)\\
    &=\frac1{\Thetaq{-q^2z}}\left(
    \frac{\frac{\lambda}{qz}\dPq{\frac{\lambda}{qz}}}{\Pq{\frac{\lambda}{qz}}}
    -\frac{\lambda\dPq{\lambda}}{\Pq{\lambda}}
    \right)=\frac{-P_q^{(d)}(z;\lambda)}{\Thetaq{-q^2z}}
\end{split}\end{equation}
We apply the operator $f(\lambda)\mapsto \frac{-1}{\ln q}\int_{q}^1 \Pq{\frac{\lambda}{z}}
f(\lambda)\frac{\intd \lambda}{\lambda}$ to both sides, and obtain,
via \eqref{IntlamdaE},
\begin{equation}
q^2\left(v_{1,\theta}(z)-v_{1,E}(z)\right)=\frac{-P_q^{(c)}(z)}{\Thetaq{-q^2z}}.
\end{equation}
It follows that \eqref{fulldifference} holds, with
\begin{align}
    P_q^{(d)}(z;\lambda)&=\frac{\lambda\dThetaq{\lambda}}{\Thetaq{\lambda}}
    -\frac{\lambda\dThetaq{\frac{\lambda}{z}}}{z\Thetaq{\frac{\lambda}{z}}}+\frac{\ln z}{\ln q}
    =\frac{\lambda\dPq{\lambda}}{\Pq{\lambda}}-\frac{\lambda\dPq{\frac{\lambda}{z}}}{z\Pq{\frac{\lambda}{z}}},\label{Pqd1}\\
    P_q^{(c)}(z)&=\frac{-1}{\ln q}\int_{q}^1 \Pq{\frac{t}{z}}
    P_q^{(d)}(z;t)\frac{\intd t}{t}
    =\frac{1}{z\ln q}\int_q^1\frac{\Pq{t}\dPq{\frac{t}{z}}}{\Pq{\frac{t}{z}}}\intd t\label{Pqc0}\\
    &=\frac{-2\pi}{\ln q}\sum_{n=1}^\infty \left(-1\right)^n\qhat^{n^2}
    \frac{\sin\left(2\pi n\frac{\ln z}{\ln q}\right)}{\sinh\left(n\ln \qhat\right)}\nonumber.
\end{align}


\section{Cauchy-Heine Asymptotics}\label{S:CH}
\begin{figure}[h]
\centering\includegraphics[width=0.3\textwidth]{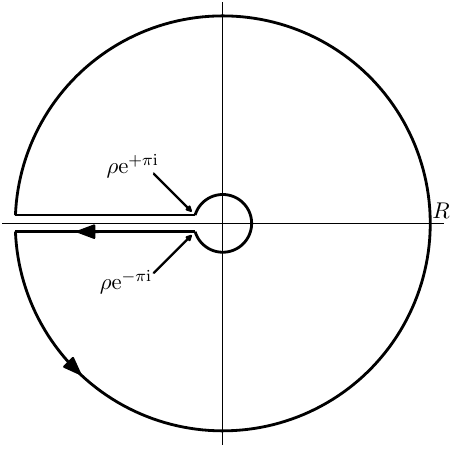}%
\caption{The Cauchy--Heine contour $\mathcal{C}$: a keyhole around the branch cut on
the negative real axis, consisting of the two line segments $\arg t=\pm\pi$ between radii
$\rho$ and $R$, closed by the circles $|t|=\rho$ and $|t|=R$. The limit $\rho\to0$ is
taken in \eqref{MB}.}

\label{fig1}
\end{figure}
Cauchy-Heine asymptotics is based on the combination of the trivial Cauchy integral representation
and an analytic continuation formula resulting in a Laplace type integral representation.

Let $0<\rho<|z|<R$, $\arg z\in(-\pi,\pi)$ and ${\mathcal C}$ the contour displayed in Figure \ref{fig1}. Then
\begin{equation}\label{MB}
\begin{split}
    y_e(z)&=\frac{1}{2\pi \iunit}\oint_{\mathcal{C}}\frac{y_e(t)}{t-z}\intd t
    =\frac{1}{2\pi\iunit}\left(
    \int_{\rho\expe^{-\pi\iunit}}^{R\expe^{-\pi\iunit}}+\int_{|t|=R}
    +\int_{R\expe^{\pi\iunit}}^{\rho\expe^{\pi\iunit}}
    +\int_{|t|=\rho}\right)\frac{y_e(t)}{t-z}\intd t\\
    &=\frac{1}{2\pi\iunit}\int_0^{R\expe^{-\pi\iunit}}\frac{y_e(t)-y_e(t\expe^{2\pi \iunit})}{t-z}\intd t
    +\frac{1}{2\pi\iunit}\int_{|t|=R}\frac{y_e(t)}{t-z}\intd t\\
   &=K(q)\int_0^R\frac{\K{q^2t}u(-t)}{z+t}\intd t
   +\frac{1}{2\pi\iunit}\int_{|t|=R}\frac{y_e(t)}{t-z}\intd t,
\end{split}
\end{equation}
where we took the limit $\rho\to 0$ and used connection relation \eqref{connection}.
In the penultimate integral, the $\K{q^2t}$ is $q$-exponentially small and we may
use the truncated geometric progression
\begin{equation}
    \frac{1}{z+t}=\sum_{n=0}^{N-1}\frac{\left(-z\right)^n}{t^{n+1}}
    +\frac{\left(-z/t\right)^N}{z+t},
\end{equation}
to obtain
\begin{equation} \label{truncated}
    y_e(z)=\sum_{n=0}^{N-1} c_n q^{-\binom{n}{2}}\left(-qz\right)^n+R_{N,e}(z),
\end{equation}
with
\begin{equation} \label{RNz}
    R_{N,e}(z)=\left(-z\right)^NK(q)\int_0^R\frac{\K{q^2t}u(-t)t^{-N}}{z+t}\intd t
   +\frac{z^N}{2\pi\iunit}\int_{|t|=R}\frac{y_e(t)t^{-N}}{t-z}\intd t,
\end{equation}
and an integral representation for the coefficients,
\begin{equation}\label{cnCH}
    c_nq^{n-\binom{n}{2}}=
    K(q)\int_0^R\K{q^2t}u(-t)t^{-n-1}\intd t
   +\frac{1}{2\pi\iunit}\int_{|t|=R}y_e(t)\left(-t\right)^{-n-1}\intd t.
\end{equation}
These representations are mainly useful in the case of large $N$ and $n$.
For example, in the case that $n$ is large, the final integral in \eqref{cnCH} is negligible
compared to the first. In the first integral in \eqref{cnCH}, the main contribution comes from the local
behaviour near $t=0$, and we obtain an asymptotic expansion when we replace $u(-t)$ by its
Maclaurin series \eqref{uexpansion} via
\begin{equation}
\begin{split}
    c_nq^{n-\binom{n}{2}}
    &\sim K(q)\sum_{m=0}^\infty \left(-1\right)^m d_m
    \int_0^\infty \K{q^2t}t^{m-n-1}\intd t\\
    &\sim K(q)\sum_{m=0}^\infty \left(-1\right)^m d_m
    q^{-\frac{1}{2}(n-m)(n-m-3)}\\
    &\sim K(q)q^{n-\binom{n}{2}}
    \sum_{m=0}^\infty \left(-1\right)^m d_m 
    q^{(n-2)m-\binom{m}{2}},
\end{split}
\end{equation}
as $n\to\infty$. We therefore have the asymptotics of the $c_n$ coefficients as
\begin{equation}\label{cnasymptotics}
    c_n\sim K(q)\sum_{m=0}^\infty \left(-1\right)^m d_m 
    q^{-\frac12(m+3)m} q^{mn}.
\end{equation}
The optimal number of terms in \eqref{truncated} is $N=N_o\sim\frac{\ln|z|}{\ln q}$. 
Hence, when $z\to 0$, $N_o$ will be large and in a similar manner we obtain
\begin{equation} \label{RNzLevel1}
\begin{split}
     R_{N_o,e}(z)&\sim \left(-z\right)^{N_o}K(q)
    \sum_{m=0}^\infty \left(-1\right)^m d_m\int_0^\infty
    \frac{\K{q^2 t}t^{m-N_o}}{z+t}\intd t\\
    &=\left(-1\right)^{N_o}K(q)
    \sum_{m=0}^\infty \left(-z\right)^m d_m
    q^{-\frac12(N_o-m)(N_o-m-5)}v_{N_o-m,e}(z),
\end{split}\end{equation}
a divergent asymptotic expansion in terms of hyperterminant functions.
We focus on the initial terms in this expansion and consider $N_o$ large and $m$ small.
In the second line of \eqref{RNzLevel1}, the factor $q^{(m+\frac52)N_o}$ will give us the
asymptotic property of this expansion.

\emph{Level $0$:} By taking only the first term in \eqref{RNzLevel1}, we obtain that the optimal remainder
of the original asymptotic expansion is
$\bigO\left(z^{\frac32}\exp\left(\frac{\ln^2|z|}{2\ln q}\right)\right)$ as $z\to 0$ and 
$N=N_o\sim\frac{\ln|z|}{\ln q}$. 

\emph{Level $1$:} We combine the truncated original expansion with re-expansion \eqref{RNzLevel1},
and take $N_o\sim\beta\frac{\ln|z|}{\ln q}$ with $\beta>1$. 
From \eqref{vnrecurrence}, we obtain that $v_{N_o-m,e}(z)\sim q^{-2}z^{N_o-m-1}$.
When we consider the $N_o$th term in re-expansion \eqref{RNzLevel1}, that is, take $m=N_o$, 
we see that it is of order $d_{N_o}z^{N_o}$. We expect that the remainder after $N_o$ terms
is of a similar size. However, $d_{N_o}z^{N_o}$ is the $N_o$th term in the convergent Taylor series
of $u(z)$. Hence, $d_{N_o}z^{N_o}\to0$ as $N_o\to\infty$. In \S\ref{S:Numerics}, we will
show that it is really possible to obtain arbitrary precision, but that it is not numerically
stable.


\section{Basic Hypergeometric Series $_2\phi_0$ and $_0\psi_2$}\label{S:2phi0}

A formal solution to the $q$-difference equation,
\begin{equation} \label{2phi0diffeq}
    xw(xq^{-2})+(q-(a+b)x)w(xq^{-1})-(q-abx)w(x)=0,
\end{equation}
is the $\qphi{2}{0}{a,b}{-}{q}{x}$ basic hypergeometric series given by
\begin{equation}
    w(x)=\qphi{2}{0}{a,b}{-}{q}{x}=\sum_{n=0}^\infty \frac{\qpr{a}{q}{n}\qpr{b}{q}{n}}{\qpr{q}{q}{n}}q^{-\binom{n}{2}}\left(-x\right)^n.
\end{equation}
Taking $z=x/q$ and $y(z)=w(x)$ we obtain the standard form
\begin{equation}
    zy(zq^{-2})+(1-(a+b)z)y(xq^{-1})-(1-abz)y(z)=0.
\end{equation}
Calculation of the $c_n$ coefficients returns the original basic hypergeometric 
series with $c_n=\frac{\qpr{a}{q}{n}\qpr{b}{q}{n}}{\qpr{q}{q}{n}}$. 
It follows from \eqref{cnasymptotics} that 
$K(q)=\frac{\qpr{a}{q}{\infty}\qpr{b}{q}{\infty}}{\qpr{q}{q}{\infty}}$.

The second solution of \eqref{2phi0diffeq} near $x=0$ is given by
\begin{equation}
\begin{split}
    w_2(x)&=\K{\expe^{\pi\iunit}qx}\qpr{abx}{q}{\infty}
    \qphi{2}{1}{\frac{q}{a},\frac{q}{b}}{0}{q}{abx}\\
    &=\K{\expe^{\pi\iunit}qx}\qpr{bqx}{q}{\infty}
    \qphi{1}{1}{\frac{q}{a}}{bqx}{q}{aqx}.
\end{split}\end{equation}
Two independent solutions near $x=\infty$ are
\begin{equation}
\begin{split}
    w_3(x)&=x^{-\ln(a)/\ln(q)}\qphi{2}{1}{a,0}{\frac{aq}{b}}{q}{\frac{q}{abx}},\\
    w_4(x)&=x^{-\ln(b)/\ln(q)}\qphi{2}{1}{b,0}{\frac{bq}{a}}{q}{\frac{q}{abx}}.
\end{split}\end{equation}
In \cite{Morita2013}, the connection relations between the solutions of \eqref{2phi0diffeq} are
discussed, focusing on the relation between $w_2(x)$, $w_3(x)$ and $w_4(x)$.

The formal $_2\phi_0$ series will have three $q$-Borel-Laplace transforms denoted by
$\UqK{\ee}{a}{b}{x}$, in which $\ee$ is either $E$, $\theta$ or $\lambda$.
More explicitly,
\begin{align}
    \UqK{\ee}{a}{b}{x}&= \int_0^{\infty}
    \qphi{2}{1}{a,b}{0}{q}{-t}\K{\frac{t}{x}}\frac{\intd t}{t},\qquad e=E,\theta,\label{2phi0BorelLaplaceK}\\
    \UqK{\lambda}{a}{b}{x}&=\sum_{n=-\infty}^\infty 
    \frac{\qphi{2}{1}{a,b}{0}{q}{-q^n\lambda}}{\Thetaq{\frac{q^n\lambda}{x}}},\label{2phi0BorelLaplaceL}
\end{align}
with the notation in \eqref{X}.
We do use the ${\rm U}$ notation to indicate that this
function is a $q$-version of the Kummer $U$ function which has the representation
$U(a,b;x)=x^{-a}{}_2F_0\left(\mytop{a,a-b+1}{-};-x^{-1}\right)$, see
\cite[\href{http://dlmf.nist.gov/13.6.E21}{13.6.21}]{NIST:DLMF}.
Hence, the $\UqK{\ee}{a}{b}{x}$ should represent an important class of special functions.

Solutions $w_3(x)$ and $w_4(x)$ show us that the growth near $x=\infty$ is algebraic.
It follows that, for the resummations of $_2\phi_0$, the final integral in \eqref{MB} vanishes
as $R\to\infty$ as long as we assume that $|a|$, $|b|<1$. Thus
\begin{equation} \label{UqCauchyHeine}
    \UqK{\ee}{a}{b}{x}=\frac{\qpr{a}{q}{\infty}\qpr{b}{q}{\infty}}{\qpr{q}{q}{\infty}}
    \int_0^\infty\frac{\qpr{-bqt}{q}{\infty}\qphi{1}{1}{\frac{q}{a}}{-bqt}{q}{-aqt}}{t+x}
    \K{qt}\intd t,
\end{equation}
when $|a|,|b|<1$. This Cauchy-Heine integral representation can be useful as an alternative
method to obtain error-bounds.

We will have a better look at the simplest case, $a=b=0$. In that case, 
$\qphi{2}{1}{a,b}{0}{q}{-t}=\frac1{\qpr{-t}{q}{\infty}}$. Thus
\begin{equation}\label{tempUq}
\begin{split}
    \UqK{\lambda}{0}{0}{x}&=\sum_{n=-\infty}^\infty 
    \frac1{\qpr{-q^n\lambda}{q}{\infty}\Thetaq{\frac{q^n\lambda}{x}}}
    =\sum_{n=-\infty}^\infty \frac{\qpr{-\lambda}{q}{n}q^{\binom{n}{2}}\left(
    \frac\lambda{x}\right)^n}{\qpr{-\lambda}{q}{\infty}\Thetaq{\frac{\lambda}{x}}}\\
    &=\sum_{n=-\infty}^\infty \frac{\qpr{-\lambda}{q}{-n}q^{\binom{n}{2}}\left(
    \frac{qx}\lambda\right)^n}{\qpr{-\lambda}{q}{\infty}\Thetaq{\frac{\lambda}{x}}}
    =\sum_{n=-\infty}^\infty \frac{q^{2\binom{n}{2}}\left(\frac{q^2x}{\lambda^2}\right)^n}{%
    \qpr{\frac{-q}\lambda}{q}{n}\qpr{-\lambda}{q}{\infty}\Thetaq{\frac{\lambda}{x}}},
\end{split}
\end{equation}
in which we have used once $n\to-n$ and the identity 
$\qpr{a}{q}{-n}=\left(-q/a\right)^nq^{\binom{n}{2}}/\qpr{q/a}{q}{n}$. We express the final
bilateral series in terms of the basic hypergeometric function 
\begin{equation}\label{psi02}
    \qpsi{0}{2}{-}{a,b}{q}{z}=\sum_{n=-\infty}^\infty
    \frac{q^{2\binom{n}{2}}z^n}{\qpr{a}{q}{n}\qpr{b}{q}{n}},
\end{equation}
and use the fact
that $\UqK{\lambda}{0}{0}{x}$ is a resummation of $\qphi{2}{0}{0,0}{-}{q}{x}$:
\begin{equation}
    \UqK{\lambda}{0}{0}{x}
    =\frac{\qpsi{0}{2}{-}{\frac{-q}\lambda,0}{q}{\frac{q^2x}{\lambda^2}}}{%
    \qpr{-\lambda}{q}{\infty}\Thetaq{\frac{\lambda}{x}}}\sim \sum_{n=0}^\infty
    \frac{q^{-\binom{n}{2}}\left(-x\right)^n}{\qpr{q}{q}{n}},
\end{equation}
as $x\to 0$. We convert this into
\begin{equation}
    \qpsi{0}{2}{-}{a,0}{q}{z}\sim\qpr{\frac{q}{a}}{q}{\infty}\Thetaq{\frac{-z}{a}}\sum_{n=0}^\infty
    \frac{q^{-\binom{n}{2}}\left(\frac{-z}{a^2}\right)^n}{\qpr{q}{q}{n}},
\end{equation}
as $z\to 0$, as long as $a\not=0$. Hence, the $z\to0$ asymptotics of this bilateral basic
hypergeometric function is a simple theta function times a divergent series.

Another interesting representation can be obtained from the final expression in \eqref{tempUq} 
by taking $\lambda=\lambda(s):=-q^{-s}$ and $x=z/q$:
\begin{equation}\begin{split}
    \UqK{\lambda(s)}{0}{0}{\frac{z}{q}}
    &=\sum_{n=-\infty}^\infty \frac{q^{n^2+2sn}z^n}{%
    \qpr{q^{1+s}}{q}{n}\qpr{q^{-s}}{q}{\infty}\Thetaq{-zq^s}}\\
    &=\frac{p(z;s)}{\Thetaq{-z}}\sum_{n=-\infty}^\infty \frac{q^{\left(n+s\right)^2}z^{n+s}}{%
    \qpr{q}{q}{n+s}}=\frac{p(z;s)}{\Thetaq{-z}}F_q(z\expe^{-\pi\iunit};s),
\end{split}\end{equation}
in which
\begin{equation}
    p(z;s)=\frac{q^{-s^2}z^{-s}\qpr{q}{q}{\infty}^2\Thetaq{-z}}{\Thetaq{-q^{-s}}\Thetaq{-zq^s}},
    \qquad F_q(z;s)=\expe^{\pi\iunit s}\sum_{n=-\infty}^\infty 
    \frac{\left(-1\right)^nq^{\left(n+s\right)^2}z^{n+s}}{\qpr{q}{q}{n+s}}.
\end{equation}
The multiplier $p(z;s)$ is $q$-periodic in $z$ and 1-periodic in $s$, and in \cite[p.~346]{Ismail2013}, it was observed
that $F_q(z;s)$ is a family of solutions of $y(z)-y(qz)+qzy\left(q^2z\right)=0$.
With the results in our paper, it is not surprising that there is a free parameter.
When we translate \eqref{fulldifference} to this specific case, we obtain the connection
relation
\begin{equation}
    p\left(z\expe^{\pi\iunit};s
    \right)F_q(z;s)=
    \Thetaq{z}\UqK{\theta}{0}{0}{\frac{z\expe^{\pi\iunit}}{q}}-
    \frac{P_q^{(d)}(z\expe^{\pi\iunit};\lambda(s))}{\qpr{q}{q}{\infty}}F_q(z;0).
\end{equation}
Note that $F_q(z;0)=\sum\limits_{n=0}^\infty \frac{q^{n^2}\left(-z\right)^n}{\qpr{q}{q}{n}}=A_q(z)$
is the Ramanujan function \cite{Ramanujan1988}.


\section{A numerical example}\label{S:Numerics}
For our numerical example, we take the $q$-difference equation
\begin{equation}
    zy(zq^{-2})+\left(1-z\right)^{\frac12}y(zq^{-1})-\left(1-\tfrac{z}2\right)^{\frac13}y(z)=0.
\end{equation}
We take this example to illustrate that the coefficients do not have to be rational functions
of $z$. The numerics will be for $q=\frac12$, $z=\frac1{1000}$ and $\lambda=\frac65$.  
We will display 8 digits. In this case, we have
$\qhat=4.2885452\times10^{-13}$ and, from \eqref{Kq}, we have $K(q)=1.5154975$.

In this case, the equation for $u(z)$ is
\begin{equation}\label{exampleUeq}
    u(zq^{-2})-\left(1-z\right)^{\frac12}u(zq^{-1})-qz\left(1-\tfrac{z}2\right)^{\frac13}u(z)=0.
\end{equation}
The function $u(z)$ is analytic near $z=0$. When we approach $z=1$ in \eqref{exampleUeq}, we obtain
from the first 2 terms that $u(zq^{-2})$ must have a square-root singularity at $z=1$, that is,
$u(z)$ has a square-root singularity at $z=q^{-2}=4$. 
This is the singularity nearest to the origin. It follows from Darboux's method that 
for the Taylor-coefficients in \eqref{uexpansion}, we have $d_n=\bigO\left(q^{2n}n^{-\frac32}\right)$
as $n\to\infty$.

Below we do sum the exact resummations. Hence, we have to evaluate
\begin{equation}
    v_{0,e}(z)=\int_0^\infty \frac{t^{-1}\K{\frac{t}{qz}}}{1+t}\intd t
    =\int_{-\infty}^\infty \frac{\K{\frac{\exp\tau}{qz}}}{1+\exp\tau}\intd \tau.
\end{equation}
We evaluated the $\tau$-integral via the trapezoidal rule, because it converges exponentially fast.

We evaluate the three resummations via the exact resummations \eqref{exactresummation} giving us
$y_e(z)=0.99933407$, with the differences $y_\theta(z)-y_E(z)=3.3647736\times 10^{-41}$
and $y_\theta(z)-y_\lambda(z)=-2.0195507\times 10^{-30}$. These differences are surprisingly small.
This follows from our formulas for the differences \eqref{fulldifference} in which
we have for the universal $q$-periodic functions
$P_q^{(c)}(z)=-7.1131070\times 10^{-25}$ and $P_q^{(d)}(z;\lambda)=4.2693156\times 10^{-14}$.
As mentioned before, these two values are small because of the smallness of $\qhat$.

\begin{figure}[h]
\centering\includegraphics[width=0.6\textwidth]{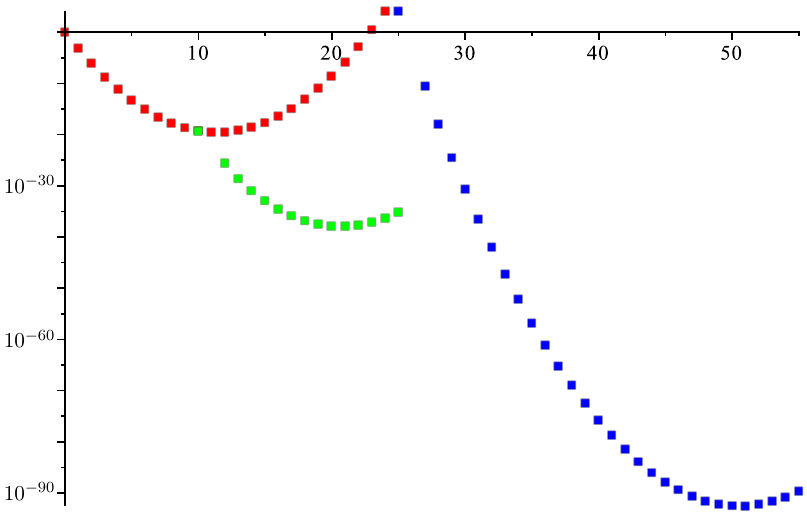}%
\caption{Magnitudes of the terms in the combined level-$1$ approximation
\eqref{Remainder}+\eqref{RNzLevel1}, on a $\log_{10}$ vertical scale against term
index. Red: terms of the original expansion \eqref{originaly}. Green: re-expansion
\eqref{RNzLevel1} with $N_o=10$ (reaching $38$ digits). Blue: re-expansion with
$N_o=25$ (reaching $92$ digits).}

\label{fig2}
\end{figure}

We will also show how these differences compare with optimal truncated asymptotics.
The optimal number of terms in the original asymptotic expansion \eqref{originaly} is
$N_o=10\approx \frac{\ln z}{\ln q}$. The remainder after taking $10$ terms is
$R_{N_o,E}(z)=3.6884443\times10^{-20}$. Hence, much larger than the differences between the three
resummations. To obtain a better asymptotic approximation, we have to re-expand the remainder
in a level 1 hyperasymptotic approximation via \eqref{RNzLevel1}. In Figure \ref{fig2}, we display
the absolute values of the terms in the level 1 hyperasymptotic approximation.
Note that, approximately, the size of the terms agree with the size of the remainders.
If we take $N_o=10$, then the red dots are the terms in the original expansion. It should be
obvious that $N_o=10$ is approximately the optimal number of terms and we obtain 20 digits
precision. The terms of the re-expansion \eqref{RNzLevel1} are the green dots and it follows that
we can obtain 38 correct digits if we take 10 terms in the re-expansion. However, it seems that
we can do much better by taking more terms in the original expansion. If we take $N_o=25$ then the
blue dots are the terms in the re-expansion with $N_o=25$. If we also take 25 terms in the
re-expansion then we can achieve 92 correct digits in our asymptotic approximation. 
We note that, when we take even more terms in the original expansion, the additional later terms will 
be larger than $\bigO(1)$ and this leads to numerical instability.

We did observe in \S\ref{S:exact} that the smallest term in the re-expansion is of order
$d_{N_o}z^{N_o}$ and, above, we observed that $d_n=\bigO\left(q^{2n}n^{-\frac32}\right)$.
Hence, the smallest term should be of size $q^{2N_o}N_o^{-\frac32}z^{N_o}=7.1\times10^{-93}$.

\appendix

\section{Theta function identities}\label{Apptheta}
In the proofs of the main theorems, we will use the following identities.
\begin{lemma}\label{lemmaparfrac}
For $x,y,x/z,x/y\notin \left\{0,-1,-q^{\pm1},-q^{\pm2},\ldots \right\}$, we have
\begin{align} 
    \frac{1}{\Thetaq{x}}&=\frac{1}{\qpr{q}{q}{\infty}^3}
    \sum_{n=-\infty}^\infty 
    \frac{\left(-1\right)^nq^{\binom{n}{2}}}{x+q^{-n}},\label{PartFrac1}\\
    \sum_{n=-\infty}^\infty \frac{1}{\left(x+q^{-n}\right)\Thetaq{yq^n}}&=
    \sum_{n=-\infty}^\infty \frac{1}{\left(y+q^{-n}\right)\Thetaq{xq^n}},\label{PartFrac1a}\\
    \frac{\Thetaq{ax}}{\Thetaq{x}}&=
    \frac{\Thetaq{-a}}{\qpr{q}{q}{\infty}^3}
    \sum_{n=-\infty}^\infty 
    \frac{a^n}{1+x q^{n}},\qquad\qquad q<|a|<1,\label{PartFrac2}\\
    \frac{\Thetaq{ax}}{\Thetaq{x}}&=
    \frac{-\pi\Thetaq{-a}}{\qpr{q}{q}{\infty}^3\ln q}
   \sum_{n=-\infty}^\infty 
    \frac{x^{-(2n\pi\iunit+\ln a)/\ln q}}{\sin\left(\frac\pi{\ln q}(2n\pi\iunit+\ln a)\right)},\label{PartFrac2a}\\
    \frac{x\dThetaq{x}}{\Thetaq{x}}-\frac{y\dThetaq{y}}{\Thetaq{y}}
    &=\sum_{n=-\infty}^\infty \left(\frac{x}{x+q^n}-\frac{y}{y+q^n}\right),\label{lemmaparfrac1}\\
    \left(\frac{x\dThetaq{x}}{\Thetaq{x}}\right)'
    &=\sum_{n=-\infty}^\infty \frac{q^n}{\left(x+q^n\right)^2},\label{lemmaparfrac2}
\end{align}
\begin{equation}\label{lemmaparfrac3}
\begin{split}
    &\frac{x\dThetaq{x}}{\Thetaq{x}}
    -\frac{x\dThetaq{x/z}}{z\Thetaq{x/z}}
    -\frac{y\dThetaq{y}}{\Thetaq{y}}
    +\frac{y\dThetaq{y/z}}{z\Thetaq{y/z}}\\
    &\qquad\qquad\qquad\qquad\qquad\qquad=\frac{\qpr{q}{q}{\infty}^3\Thetaq{-z}\Thetaq{\frac{-x}{y}}
    \Thetaq{\frac{-xy}{z}}}{\Thetaq{x}\Thetaq{y}
    \Thetaq{\frac{x}{z}}\Thetaq{\frac{z}{y}}}.
\end{split}\end{equation}
\end{lemma}
\begin{proof}
Expansion \eqref{PartFrac1} can be obtained by considering the integral 
$\int_0^\infty\frac{\intd z}{(z-x)\Thetaq{z}}$ and rotating the contour of integration by $2\pi$.
For \eqref{PartFrac1a}, use \eqref{PartFrac1} on both sides.
Equation \eqref{PartFrac2} follows from \cite[Eq.(10.6.1)]{AAR1999}.
For \eqref{PartFrac2a}, we start with the constraints $x,a\in(q,1)$ and consider the residues of the
Mellin--Barnes type integral
\begin{equation}\label{tempMB}
\begin{split}
    &\frac{-\pi}{2\pi\iunit}\left(\int_{M+\frac12-\iunit\infty}^{M+\frac12+\iunit\infty}-
    \int_{-N-\frac12-\iunit\infty}^{-N-\frac12+\iunit\infty}\right)
    \frac{x^\tau}{\left(1-aq^\tau\right)\sin\pi\tau}\intd\tau\\
    &\qquad\qquad=
    \frac{-\pi}{\ln q}\sum_{n=-\infty}^\infty 
    \frac{x^{-(2n\pi\iunit+\ln a)/\ln q}}{\sin\left(\frac\pi{\ln q}(2n\pi\iunit+\ln a)\right)}
    -\sum_{m=-N}^M \frac{\left(-x\right)^m}{1-aq^m}.
\end{split}\end{equation}
We push $M,N\to\infty$. The left-hand side of \eqref{tempMB} vanishes, and the final term of \eqref{tempMB}
can be expressed in the required quotient of theta functions via \eqref{PartFrac2}.

Let $N$ be a positive integer and define $\Thetaq{\tau;N}=\qpr{q,-\tau,-q/\tau}{q}{N}$. 
The reader can check that
\begin{equation*}
    \frac{x\dThetaq{x;N}}{\Thetaq{x;N}}-\frac{y\dThetaq{y;N}}{\Thetaq{y;N}}
    =\sum_{n=1-N}^N \left(\frac{x}{x+q^n}-\frac{y}{y+q^n}\right).
\end{equation*}
Now let $N\to\infty$. For \eqref{lemmaparfrac2}, consider $\lim_{y\to x}\frac{\eqref{lemmaparfrac1}}{x-y}$.

We take $z=\expe^{2\iunit \zeta}$. 
As functions of $\zeta$, both expressions of \eqref{lemmaparfrac3} are elliptic
functions of order 2 with periods $\pi$ and $\frac{\ln q}{2\iunit}$. It is easy to show that both expressions 
of \eqref{lemmaparfrac3} have the same poles, zeros, and the same residues at the poles. 
Hence, the quotient of both sides is constant, and because of the same residues, this constant is unity.
\end{proof}

An additional identity obtained for the ratio of $\Eq{\tau}$ functions is
\begin{equation} \label{Eqab}
\frac{\Eq{at}\Eq{bt}}{\Eq{t}\Eq{abt}}=\expe^{-\ln(a)\ln(b)/\ln(q)}.
\end{equation}

For the function $\Pq{\tau}$ defined in \eqref{defPq}, we have the expansion
\begin{equation}\label{qperiodic}
\Pq{\tau}=-\ln(q)\Cq\Eq{\tau}\Thetaq{\tau}=
\qpr{\qhat^2}{\qhat^2}{\infty}
    \exp\left(\sum_{n=1}^\infty\frac{\cos\left(2\pi n\frac{\ln\tau}{\ln q}\right)}{n%
    \sinh\left(n\ln \qhat\right)}
    \right).
\end{equation}
This expansion is an alternative version of \cite[(3.7)]{Joshi2025}.
We also have
\begin{equation} \label{Pq}
    \Pq{q^t}=-\ln(q)\Cq \sum_{n=-\infty}^\infty q^{\frac12\left(t-\frac12+n\right)^2},
\end{equation}
with Fourier series,
\begin{equation} \label{PqFourier}
    \Pq{q^t}=1+2\sum_{n=1}^\infty \left(-1\right)^n\qhat^{n^2} \cos(2n\pi t)=\sum_{n=-\infty}^\infty (-1)^n\qhat^{n^2}\expe^{2n\pi\iunit t}=\theta_4(\pi t,\qhat).
\end{equation}
This clearly demonstrates that $\Pq{\tau}-1$ is extremely small for $q$ bounded away from $0$.
It follows from Fourier series \eqref{PqFourier} that
\begin{equation} \label{PqInt}
    \frac{-1}{\ln q}\int_q^1 \Pq{z\tau}\frac{\intd\tau}\tau=\int_0^1\Pq{zq^t}\intd t
    =\int_0^1\Pq{q^{t+\frac{\ln z}{\ln q}}}\intd t=\int_0^1\Pq{q^t}\intd t=1.
\end{equation}
We may also derive, from \cite[Example 14, p.\ 489]{WW1927},
\begin{equation}\label{PqReciprocal}
\frac1{\Pq{q^t}}=\frac1{\qpr{\qhat^2}{\qhat^2}{\infty}^3}
\left(\widetilde{a}_0+2\sum_{n=1}^\infty \widetilde{a}_n\qhat^n\cos(2n\pi t)\right),
\end{equation}
with
\begin{equation}\label{tan}
    \widetilde{a}_n=\sum_{m=0}^\infty \left(-1\right)^m
    \qhat^{m(m+2n+1)}
    =\qphi{1}{1}{\qhat^2}{0}{\qhat^2}{\qhat^{2n+2}}.
\end{equation}
These coefficients have the following recurrence and normalising relations,
\begin{equation}
    \qhat^{2n+2}\widetilde{a}_{n+1}+\widetilde{a}_n=1,\qquad 
    \sum_{n=0}^\infty \frac{\left(-1\right)^n\qhat^{n(n+1)}\widetilde{a}_n}{
    \qpr{\qhat^2}{\qhat^2}{n}}=\qpr{\qhat^2}{\qhat^2}{\infty}^2.
\end{equation}
The function $\Eq{\tau}$ and hence $\Pq{\tau}$ are multi-valued:
\begin{equation}
    \frac{\Eq{\tau\expe^{2\pi\iunit}}}{\Eq{\tau}}=\frac{\Pq{\tau\expe^{2\pi\iunit}}}{\Pq{\tau}}=
    -\qhat^{-1}\tau^{\frac{2\pi\iunit}{\ln q}}.
\end{equation}
Using the translation of the argument by a half-period, see
\cite[\href{http://dlmf.nist.gov/20.2.E10}{20.2.10}]{NIST:DLMF}, we obtain
\begin{equation}\label{PqMinus}
\Pq{\expe^{\pi\iunit}q^t}
=-\iunit\expe^{\pi\iunit t}\qhat^{-\frac14}\theta_1(\pi t,\qhat).
\end{equation}
The function $\Pq{t}$ is the bridge between the two continuous $q$-Laplace transforms:
\begin{equation}
\begin{split}
    (\LqK{E} B)(z)=\Cq\int_0^\infty B(t) \Eq{\frac{t}{z}}\, \frac{\intd t}{t}&=
    \frac{-1}{\ln q}\int_0^\infty \frac{B(t)}{\Thetaq{\frac{t}{z}}} 
    \Pq{\frac{t}{z}}\, \frac{\intd t}{t},\\
    (\LqK{\theta} B)(z)=\frac{-1}{\ln q}\int_0^\infty \frac{B(t)}{\Thetaq{\frac{t}{z}}} 
    \, \frac{\intd t}{t}&=
    \Cq\int_0^\infty \frac{B(t) \Eq{\frac{t}{z}}}{\Pq{\frac{t}{z}}}\, \frac{\intd t}{t}.
\end{split}
\end{equation}
Finally, we observe that a quotient of two $\theta_q$ functions is a simple power times a
$q$-periodic function.
\begin{equation}\label{thetaquotient}
    \frac{\Thetaq{ax}}{\Thetaq{bx}}=x^{\frac{\ln(b/a)}{\ln q}}\frac{\Eq{b}\Pq{ax}}{\Eq{a}\Pq{bx}}.
\end{equation}

\section{Lemmas}\label{AppLemma}
\begin{lemma}\label{lemmacn}
    Assume $0<q<1$, and let $C$ be a positive constant such that $|a_n|+|b_n|\leq C^n$ for $n\geq1$.
    Then $|c_n|\leq\qpr{-2C}{q}{n}/\qpr{q}{q}{n}$ for all $n$.
\end{lemma}
\begin{proof}
    We take $e_0=1$ and 
    \begin{equation} 
    e_n=\frac{1}{1-q^n}\left(e_{n-1}+\sum_{m=0}^{n-1}e_mq^{\binom{n}{2}-\binom{m}{2}}C^{n-m}\frac{\qpr{q}{q}{m}}{\qpr{q}{q}{n-1}}\right),\qquad n\geq1.
    \end{equation}
    We compare this recurrence with \eqref{cnrecurrence}, the one for the $c_n$, and obtain that
    $|c_n|\leq e_n$ for all $n$. The formal series 
    $e(z)=\sum\limits_{n=0}^\infty e_n\qpr{q}{q}{n}q^{-\binom{n}{2}}z^n$ satisfies
    \begin{equation} 
    2e(t)=te(t/q)+\frac{e(t)}{1-Ct}+1.
    \end{equation}
    We multiply both sides by $1-Ct$ and obtain the new recurrence relation
    \begin{equation}\label{enrecurrenceNew}
    e_n=\frac{\left(1+2Cq^{n-1}\right)}{(1-q^n)}e_{n-1}-\frac{Cq^{n-1}}{(1-q^{n-1})(1-q^n)}e_{n-2},\qquad n\geq2,
    \end{equation}
    with $e_0=1$ and $e_1=(1+C)/(1-q)<\qpr{-2C}{q}{1}/\qpr{q}{q}{1}$. We obtain from \eqref{enrecurrenceNew} and induction
    that
    \begin{equation} 
    e_n\leq\frac{\left(1+2Cq^{n-1}\right)}{(1-q^n)}e_{n-1}\leq \ldots\leq \frac{\qpr{-2C}{q}{n}}{\qpr{q}{q}{n}} e_0.
    \end{equation}
\end{proof}

\begin{lemma}\label{lemmadn}
    Assume $0<q<1$, and let $A, B$ be positive constants such that $|a_m|\leq A^m$, $|b_m|\leq B^m$ 
    for all $m$. Take $C=\max(A,B)$, $r=1+qC^{-1}$ and define $e_0=1$ and 
    $e_n=\frac{rq^{2n}C^n\qpr{-rq}{q}{n-1}}{\qpr{q}{q}{n}}$ for $n\geq1$. Then $|d_n|\leq e_n$ for all $n$.
\end{lemma}
\begin{proof}
Note that $\displaystyle e_n=\frac{r\,q^{2n}}{\qpr{q}{q}{n}}\sum_{k=0}^{n-1}
q^{-k}\qpr{q}{q}{k} C^{n-k}e_k$ and that
$\left|a_{n-k}+q^{k+1}b_{n-k-1}\right|\leq C^{n-k}+qC^{n-k-1}=rC^{n-k}$.
For $n=0$, we have $d_0=e_0=1$ and therefore $|d_0|\leq e_0$.
Now assuming $|d_k|\leq e_k$ for $k=0,\ldots,n-1$, we have,
\begin{equation}
    \left|d_n\right|\leq \frac{r\,q^{2n}}{\qpr{q}{q}{n}}\sum_{k=0}^{n-1}
q^{-k}\qpr{q}{q}{k} C^{n-k}\left|d_k\right|
    \leq \frac{r\,q^{2n}}{\qpr{q}{q}{n}}\sum_{k=0}^{n-1}
q^{-k}\qpr{q}{q}{k} C^{n-k}e_k=e_n.
\end{equation}
Therefore, by strong induction, we have $|d_n|\leq e_n$ for all $n$.
\end{proof}

The function $u(z)$ will therefore have a convergent expansion on a disc with a non-zero radius as,
\begin{equation}
\begin{split}
    |u(z)|&=\left|\sum_{n=0}^\infty d_n z^n\right|
    \leq\sum_{n=0}^\infty\left|d_n\right|\left|z\right|^n
    \leq\sum_{n=0}^\infty e_n\left|z\right|^n\\
    &=1+r\sum_{n=1}^\infty\frac{\qpr{-rq}{q}{n-1}}{\qpr{q}{q}{n}}\left(Cq^2|z|\right)^n
    =\frac{1}{1+r}+\frac{r}{1+r}\sum_{n=0}^\infty\frac{\qpr{-r}{q}{n}}{\qpr{q}{q}{n}}\left(Cq^2|z|\right)^n\\
    &=\frac{1}{1+r}+\frac{r}{1+r}\frac{\qpr{-rCq^2|z|}{q}{\infty}}{\qpr{Cq^2|z|}{q}{\infty}},
\end{split}
\end{equation}
assuming $Cq^2|z|<1$.
Therefore the expansion for $u(z)$ will be convergent on disc of radius at least $|z|=\frac{1}{Cq^2}$.



\end{document}